\documentclass[10pt]{article}

\usepackage[dvips]{graphicx}
\usepackage{psfrag}
\usepackage{amsopn}
\usepackage{amstext}
\usepackage{amsfonts}
\usepackage{amssymb}
\usepackage{amscd}
\usepackage{amsthm}
\usepackage{amsfonts}
\usepackage[centertags,reqno]{amsmath}
\usepackage{latexsym}

\title{SELF--SCALED BARRIER FUNCTIONS ON SYMMETRIC CONES AND THEIR 
CLASSIFICATION}
\author{Raphael A.\ Hauser \footnote{Department of Applied Mathematics
and Theoretical Physics (DAMTP), University of Cambridge, Silver
Street, Cambridge CB3 9EW, England. E-mail: {\em
rah48@damtp.cam.ac.uk}. Research supported by the Swiss 
National Science Foundation and the Swiss Academy of Technical
Sciences, and by the Engineering and Physical Sciences Research 
Council of the UK, EPSRC grant GR/M30975.} , 
Osman G\"uler \footnote{Department of Mathematics and Statistics, 
University of Maryland Baltimore County, Baltimore, Maryland 21250, USA. 
E-mail: {\em guler@math.umbc.edu}. Research on this paper 
was conducted when the author was visiting the Industrial Engineering 
Department at Bilkent University, Ankara, Turkey. Research partially 
supported by the National Science Foundation under grant DMS--0075722.}}   

\numberwithin{equation}{section}

\DeclareMathOperator{\GL}{GL}
\DeclareMathOperator{\I}{I}
\DeclareMathOperator{\determinant}{det}
\DeclareMathOperator{\e}{e}
\DeclareMathOperator{\Int}{int}
\DeclareMathOperator{\relint}{ri}
\DeclareMathOperator{\Aut}{Aut}
\DeclareMathOperator{\Orth}{O}
\DeclareMathOperator{\Span}{span}
\DeclareMathOperator{\codim}{codim}

\theoremstyle{plain}
\newtheorem{theorem}{Theorem}[section]
\newtheorem{lemma}[theorem]{Lemma}

\newtheorem{proposition}[theorem]{Proposition}
\theoremstyle{definition}
\newtheorem{definition}[theorem]{Definition}
\newtheorem{remark}[theorem]{Remark}

\begin{document}

\thispagestyle{empty}

\vspace*{45pt}

\begin{center}
{\Large UNIVERSITY OF CAMBRIDGE}

\vspace*{6mm}
{\large Numerical Analysis Reports}

\vspace*{30mm}
\parbox[t]{105mm}{\large\bf\begin{center}
SELF--SCALED BARRIER FUNCTIONS ON SYMMETRIC CONES AND THEIR 
CLASSIFICATION
\end{center}}

\vspace*{15mm}
{\bf Raphael Hauser and Osman G\"uler}

\vspace*{45mm}
DAMTP 2001/NA03\\ 
March 2001

\vspace*{-330pt}
\begin{picture}(370,340)(0,0)
\put (3,3) {\line(1,0){364}}
\put (3,3) {\line(0,1){334}}
\put (367,3) {\line(0,1){334}}
\put (3,337) {\line(1,0){364}}
\thicklines
\put (0,0) {\line(1,0){370}}
\put (0,0) {\line(0,1){340}}
\put (370,0) {\line(0,1){340}}
\put (0,340) {\line(1,0){370}}
\end{picture}

\vfill
{\large Department of Applied Mathematics and Theoretical
Physics\\Silver Street\\Cambridge\\England CB3 9EW}
\end{center}

\thispagestyle{empty}

\thispagestyle{empty}
\newpage

\hspace{3cm}

\vfill
\thispagestyle{empty}


\newpage
\setcounter{page}{1}
\maketitle

\vspace{-1cm}

\begin{abstract} 
Self--scaled barrier functions on self--scaled cones were introduced 
through a set of axioms in 1994 by Y.~E. Nesterov and M.~J. Todd 
as a tool for the construction of long--step interior point algorithms. 
This paper provides firm foundation for these objects by exhibiting 
their symmetry properties, their intimate ties with the symmetry
groups of their domains of definition, and subsequently their
decomposition into irreducible parts and algebraic classification
theory. In a first part we recall the characterisation of the family of 
self--scaled cones as the set of symmetric cones and develop a  
primal--dual symmetric viewpoint on self-scaled barriers, results that
were first discovered by the second author. We then show in a short, 
simple proof that any pointed, convex cone decomposes into a direct sum of 
irreducible components in a unique way, a result which can also be of 
independent interest. We then show that any self--scaled barrier 
function decomposes in an essentially unique way into a direct sum of 
self--scaled barriers defined on the irreducible components of the 
underlying symmetric cone. Finally, we present a complete algebraic 
classification of self--scaled barrier functions using the 
correspondence between symmetric cones and Euclidean Jordan algebras.
\end{abstract}

\vspace{1.5in}

\noindent{\bf Abbreviated title.}  Self--scaled barrier functions.\\
\noindent{\bf AMS(MOS) subject classifications:}  
primary 90C25, 90C60, 52A41; secondary 90C06, 52A40.


\newpage

\section{Introduction}\label{sec:introduction}

In recent years a theory of interior--point methods for linear, 
semidefinite, and second--order cone programming has been developed 
within the unified framework of {\em self--scaled conic programming}. 
The origins of this theory can be traced to the works 
\cite{NT97,NT98,Guler96}. The importance of the problems which can be 
cast in this framework, and the fact that it is possible to develop 
efficient {\em long--step} interior--point methods for these problems 
have contributed to the popularity of the subject in the optimization 
community and beyond.

In order to facilitate our exposition, we consider the following pair 
of convex programs in conic duality
\begin{align}\label{problem}
\text{(P)}\quad\inf\:&\langle{x,s_0}\rangle&\text{(D)}\quad
\inf\:&\langle{x_0,s}\rangle\\
x&\in\bigl(L+x_0\bigr)\cap K&s&\in\bigl(L^\perp +s_0\bigr)
\cap K^*.\nonumber
\end{align}
Here $E$ is a  finite dimensional Euclidean space equipped with an
inner product $\langle{\cdot,\cdot}\rangle$, $L$ is a linear subspace 
of $E$, and $L^\perp$ its 
orthogonal complement. The cone $K$ is a {\em regular} 
(closed, convex, pointed, solid) cone, $x_0\in\Int(K)$, and 
$s_0\in\Int(K^*)$, where $K^*$ is the dual cone 
\begin{equation}\label{eq:dualcone}
K^* := \{s\in E: \langle{x,s}\rangle\ge0,\quad\forall x\in K\}.
\end{equation}

Interior--point algorithms can be used to solve \eqref{problem} over {\em any} 
regular cone, provided one has a {\em self--concordant} barrier 
function $F(x)$ defined over the interior $\Int(K)$ of $K$. The reader
may consult the 
authoritative monographs \cite{NN94,Renegar99} for a detailed 
treatment of self--concordant functions and interior--point methods. 
In a generic self--concordant barrier function, one has control over 
the behaviour of the Hessians $F''(y)$ only when $y$ lies in the 
local ball $\bigl\{y:\langle F''(x)(y-x),y-x\rangle<1\bigr\}$, 
leading to {\em short--step} interior--point methods. Although these 
methods have a polynomial running--time guarantee, they tend to be 
less efficient linear programming solvers in practice than 
{\em long--step} interior--point methods. The theoretical basis 
for this latter type of algorithm is the fact that the self--concordant 
barrier function \[x\mapsto -\sum_{i=1}^n\ln x_i\] has additional
properties 
which make it possible to control it in {\em all} of $\Int(K)$. 

In \cite{NT97}, Nesterov and Todd isolate two properties of the 
barrier $-\sum_{i=1}^n\ln x_i$ which they identify as being
responsible for making the long--step approach succeed in linear 
programming. They generalise these properties (see \eqref{eq:ss-1} 
and \eqref{eq:ss-2}) and call self--concordant barrier functions 
satisfying these conditions {\em self--scaled}. Since these properties 
also impose certain conditions on the domain of definition of such
functions, Nesterov and Todd call the closures of such domains 
{\em self--scaled cones}. For convenience, we recall these concepts 
here: 

\begin{definition}
Let $K\subseteq E$ be a regular cone. A self--concordant barrier
function $F:\Int(K)\to{\mathbb R}$ is called {\em self--scaled} if
$F''(x)$ is non--singular for every $x\in K$, $F$ is logarithmically
homogeneous, that is, there exists a constant $\nu>0$ such that 
\begin{equation}
F(tx) = F(x)-\nu\ln t,\qquad\forall x\in\Int(K),\ t>0,
\label{eq:log-homogeneous}
\end{equation}
and if $F$ satisfies the following two properties 
\begin{align}
F''(w)x&\in\Int(K^*),\qquad\forall x,w\in\Int(K),\label{eq:ss-1} \\
F_*\bigl(F''(w)x\bigr)&=F(x)-2F(w)-\nu,\qquad\forall x,w\in\Int(K).
\label{eq:ss-2}
\end{align}
If $K$ allows such a barrier function, then $K$ is called a {\em
self--scaled cone}. 
\end{definition}
The {\em dual barrier} $F_*:\Int(K)\to{\mathbb R}$ that appears in 
the last Axiom \eqref{eq:ss-2} is defined as 
$F_*(s):=\sup\bigl\{-\langle{x,s}\rangle-F(x):x\in\Int(K)\bigr\}$. 
Theorem 3.1 in \cite{NT97} states that \eqref{eq:ss-1} can be 
strengthened to 
\begin{theorem}
\label{thm:scaling}
If $x\in\Int(K)$ and $y\in\Int(K^*)$, then there exists a 
unique point $w\in\Int(K)$ such that \[F''(w)x=y.\] 
Moreover, if $w\in\Int(K)$ then $F''(w)(K) = K^*$.\qed
\end{theorem}

The point $w$ is called the {\em scaling point} of $x$ and $y$. 
The last statement is a consequence of the first part of the theorem
and of Equation (3.2) from Nesterov and Todd's paper \cite{NT97}. 
We reproduce this formula here for convenience:
\begin{equation}\label{eq:ff1stversion}
F''(x)=F''(w)F_*''\bigl(F''(w)x\bigr)F''(w).
\end{equation}
See also Lemma~\ref{lemma:fund-eq} below, where this formula reappears
and where we give a proof of this important identity.
 
We would like to mention that Rothaus \cite{Rothaus60}, using rather 
elementary tools, proves a number of results which are useful in 
Section~\ref{sec:group} of this paper. 
Two {\em key} results are \cite{Rothaus60}, 
Theorem 3.12 and Corollary 3.15, p.\ 205. These results imply 
Theorem~\ref{thm:scaling} for the universal barrier function, a
special self--concordant barrier function defined by Nesterov and
Nemirovskii 
\cite{NN94} which is further discussed below. Theorem~\ref{thm:scaling} 
is an independently discovered extension of Rothaus's result to all 
self--scaled functions. It can be shown that, when suitably modified, 
all results of Section III in \cite{Rothaus60} can be extended to
general self--scaled barriers. In particular this is true for 
Theorem~\ref{thm:scaling}. 

Nesterov and Todd~\cite{NT97,NT98} demonstrate that self--scaled 
barrier functions can indeed be used to develop various long--step 
interior--point methods for linear optimization over self--scaled 
cones, in particular for semidefinite programming and for convex 
quadratic programming with convex quadratic constraints. 

Inspired by the paper of Vinberg~\cite{Vinberg63}, O.\ 
G\"uler~\cite{Guler96} develops the relationship between the 
{\em universal barrier function} of Nesterov and Nemirovskii 
\cite{NN94} and the {\em characteristic function} of the cone $K$, 
\begin{equation}\label{charfn}
\varphi_K(x):=\int_{K^*}\e^{-\langle x,y\rangle}dy, 
\end{equation}
introduced in 1957 by Koecher, see \cite{Koecher99}. In particular,
G\"uler ~\cite{Guler96} shows that if $K$ is a regular cone, then 
the universal barrier function $U(x)$ for $K$ satisfies the equation   
\begin{equation}
U(x) = c_1\ln\varphi_K(x) + c_2,\label{eq:U(x)}
\end{equation} 
for some constants $c_1>0$ and $c_2$. 

Through G\"uler's paper \cite{Guler96} the concepts of 
{\em homogeneous cones}, {\em homogeneous self--dual cones} (or 
{\em symmetric cones}), {\em Euclidean Jordan algebras} and 
{\em Siegel domains} as well as the classification theory of 
symmetric cones and Euclidean Jordan algebras, known to mathematicians
since 1960 and 1934 respectively, were first introduced into the in the 
interior--point literature. The interested 
reader can find a complete treatment of these classification results 
in the book of J.\ Faraut and A.\ Kor\'anyi~\cite{FK94}. See also
\cite{Koecher99} for a different treatment of 
some of the same topics. Because of their importance for this paper, 
we recall some of the concepts mentioned above. 
\begin{definition} 
Let $K\subseteq E$ be a regular cone. The {\em automorphism group} of $K$ 
is the set of all non--singular linear maps $A:E\to E$ that leave $K$ 
invariant, i.e.,
\begin{equation*}
\Aut(K):=\bigl\{A\in\GL(E):A(K)=K\bigr\} 
\end{equation*}
The cone $K$ is called {\em homogeneous} if $\Aut(K)$ acts
transitively on $K$, 
that is, given arbitrary points $x,y\in K$, there exists a map
$A\in\Aut(K)$ 
such that $Ax=y$. The cone $K$ is called {\em self--dual} if $E$ can be 
endowed with an inner product such that $K^*=K$ where $K^*$ is defined
with 
respect to this inner product, see \eqref{eq:dualcone}. The cone $K$
is called {\em homogeneous self--dual} if $K$ is both homogeneous and 
self--dual. 
\end{definition}
Homogeneous self--dual cones are also called {\em symmetric cones} 
in \cite{FK94}, a terminology which we shall adopt in this paper. 
 
The motivation behind \cite{NT97,NT98} in contrast to
\cite{Guler96} is rather different: While the first two papers deal 
with long--step interior--point 
methods and regular cones on which such methods can be designed, the
latter one deals with the universal barrier function and the symmetry 
properties of regular cones, both in the group theoretic and the
duality theoretic sense. Shortly after the announcement of the paper 
\cite{NT97} G\"uler~\cite{Guler94} discovers that the families of 
self--scaled cones and of symmetric cones are identical, thus 
establishing a connection between these two previously distinct ideas. 

As mentioned earlier, symmetric cones are fully classified in the theory 
of Euclidean Jordan algebras, see \cite{Koecher99,FK94} and the references 
therein. According to this theory, each symmetric cone has a unique 
decomposition into a direct sum of elementary building blocks, so--called 
{\em irreducible} symmetric cones, of which there exist only five types. 
Three examples of symmetric cones are of particular interest to the 
optimization community: The non--negative orthant $K={\mathbb R}^n_+$,
the cone $K=\Sigma^{+}_{n}$ of $n\times n$ symmetric, positive 
semidefinite matrices over the real numbers, and the Lorentz cone 
$K=\bigl\{\bigl(\begin{smallmatrix}\tau\\x
\end{smallmatrix}\bigr)\in{\mathbb R}^{n+1}:\tau\geq\|x\|_2\bigr\}$.
The general self--scaled conic optimization problems associated with 
these cones are respectively linear programming, semidefinite
programming and second order cone programming. The latter can be seen
as a reformulation of convex quadratic programming with convex
quadratic constraints. Considering more
general symmetric cones, one can treat linear optimization problems
with mixed linear, semidefinite and convex quadratic constraints in a
single unified framework. 

Motivated by \cite{Guler96} and by the fact that only 
a small number of examples of self-scaled barrier functions are
explicitly known, Hauser develops a partial algebraic 
classification theory for self--scaled barrier functions in a chapter of
his thesis ~\cite{Hauserthesis} and later announces a report 
\cite{Hauser99} based on these result. Hauser shows that any
self--scaled barrier over a symmetric cone $K$ has an essentially
unique decomposition into a direct sum of self--scaled barriers
defined on the irreducible summands of $K$. 
Using this decomposition, he classifies in particular the family of 
{\em isotropic} self--scaled barrier functions which are characterised
by rotational invariance. The insight gained from a lemma leads Hauser
to conjecture that all self--scaled barrier functions defined on 
{\em irreducible} symmetric cones must be isotropic. Hauser points out
that this conjecture, if true, would settle the 
classification problem of general self--scaled barrier functions. 
This conjecture also implies that all self--scaled barriers over 
irreducible symmetric cones are of the form 
$c_1\ln\varphi_K+c_2$, where $c_1>0$ and $c_2$ are constants. 
It follows from this theory that 
the full set of self--scaled barrier functions is readily known, and 
that all of these functions are just minor transformations of the
universal barrier function. 

In a second report \cite{Hauser00}, Hauser proves this conjecture in
the special case where $K$ is the positive semidefinite cone, i.e., he
shows that all self--scaled barrier functions for use in 
semidefinite programming are isotropic and essentially identical to 
the standard logarithmic barrier function. 
Starting from first principles, Hauser shows that the orientation
preserving part of the automorphism group of the positive semidefinite
symmetric cone is generated by the Hessians of an arbitrary
self--scaled barrier function defined on this cone, see
\cite{Hauser00}, Corollary 4.3. Hauser's solution of the isotropy 
conjecture in this special case relies on Proposition 3.3, which also 
forms the key mechanism in the proof of Corollary 4.3 in \cite{Hauser00}. 

Hauser's Corollary 4.3 is essentially a rediscovery of a result by 
Koecher, Theorem 4.9 (b), 
pp.\ 88--89 ~\cite{Koecher99}, for the special case of the cone of 
positive semidefinite symmetric matrices. 
Shortly after the announcement of Hauser's report \cite{Hauser00} 
in March of 2000, Y.\ Lim \cite{LimPC} settles the general case of 
the isotropy 
conjecture by generalising Hauser's proof while refereeing the paper. 
Subsequently, both O.\ G\"uler~\cite{Guler00} and 
S.\ Schmieta~\cite{Schmieta00} independently of each other and 
independently of Y.\ Lim prove the isotropy conjecture in the general
irreducible case. It is interesting to note that all three approaches 
to the general case of the conjecture, as well as Hauser's approach to
the special case relevant to semidefinite programming, rely on the
same deeper principle provided by Koecher's Theorem 4.9 (b) cited
above. See also Remark~\ref{rem:lim} of this article. 

The present article is a revision of Hauser's original report 
\cite{Hauser99} which is based on parts of his thesis \cite{Hauserthesis}, 
while incorporating the solution to the general problem using 
G\"uler's approach. Schmieta's report~\cite{Schmieta00} constitutes
the first document where a proof of the general classification 
Theorem~\ref{thm:main} became publicly available. 

The rest of the paper is organised as follows. In 
Section~\ref{sec:symmetries} 
we reconsider self--scaled cones and self--scaled barriers from a 
symmetric point of view. Section~\ref{sec:group} is devoted to certain
properties of self--scaled barriers which link self--scaled barriers 
to the symmetry group of their domain of definition. These results 
are needed in later sections. 
In Section~\ref{sec:decomp}, we show that any pointed, convex cone
has a unique decomposition as a direct sum of irreducible
components. This result, of which we manage to locate only 
technically more involved generalisations, may be of independent
interest. We therefore include a simple proof. We then use this
decomposition result to show that any self--scaled barrier defined on
a symmetric cone $K$ decomposes in an essentially unique way into a
direct sum of self--scaled barriers defined on the irreducible
components of $K$, which also shows that the irreducible components 
are symmetric cones themselves. This decomposition reduces the problem of
classifying self-scaled barriers to the case where the domain of
definition is irreducible, a problem we solve in
Section~\ref{sec:classification}. Theorem~\ref{thm:main} constitutes 
the main and final result of this paper. We thus present all the 
essential elements of the theory of self--scaled barrier functions in 
a single document. 

The following basic properties of 
$\nu$\nobreakdash--self--concordant logarithmically homogeneous
barrier functions and their duals will be used frequently in later 
sections. These properties are easy consequences of
logarithmic homogeneity, see \cite{NN94} or \cite{NT97}:

\begin{proposition}\label{prop2.1}
Let $\,F\,$ be a $\,\nu$-self-concordant 
logarithmically homogeneous barrier function on the regular cone 
$\,K\subset E$, and let $\,x\in\Int(K)$, $s\in\Int(K^*)$. Then 
\begin{align*}
i)&\;-F'(x)=F''(x)x\in\Int(K^*),&
ii)&\;-F_*'\bigl(-F'(x)\bigr)=x, \\
iii)&\;F_*''(-F'(x))=F''(x)^{-1},&
iv)&\;\langle x,-F'(x)\rangle=\nu,\\
v)&\;F^{(k)}(tx)=t^{-k}F^{(k)}(x),\ \forall\,t>0, k=1,2,&
vi)&\;F_*\bigl(-F'(x)\bigr)=-\nu-F(x), \\
vii)&\;F(x)+F_*(s)\ge-\nu-\nu\log\nu-\nu\log\langle{x,s}\rangle,
\end{align*}
where $F^{(1)}(x)=F'(x)$, $F^{(2)}(x)=F''(x)$ in $v)$. \qed
\end{proposition}

\section{A Symmetric View on Self--Scaledness}\label{sec:symmetries}

In this section we undertake a study of self--scaled cones and 
barrier functions while emphasising their symmetry properties in a 
duality-theoretic sense. 

Let $F$ be a self--scaled barrier function on a regular cone $K$ in a 
finite dimensional Euclidean space $E$ equipped with an inner product 
$\langle{\cdot,\cdot}\rangle$. With a given arbitrary point $e\in\Int(K)$ 
we associate an inner product 
\[\langle{u,v}\rangle_e:=\langle{F''(e)u,v}\rangle.\] 
The following result is due to G\"uler~\cite{Guler94}: 

\begin{theorem}\label{thm:K-symmetric} 
The cone $K$ is symmetric,  and $F$ is self--scaled under 
$\langle{\cdot,\cdot}\rangle_e$. 
\end{theorem}

\begin{proof}
We have 
\begin{align*} 
K_e^*&:=\bigl\{y: \langle{x,y}\rangle_e\ge0,\ \forall x\in K\bigr\} 
=\bigl\{y:\langle{F''(e)x,y}\rangle\ge0,\ \forall x\in K\bigr\}\\
&=\bigl\{y:\langle{z,y}\rangle\ge0,\ \forall z\in K^*\bigr\}=(K^*)^*=K,
\end{align*}
where the third equality follows from Theorem~\ref{thm:scaling}. Note that 
\[\langle{F''(x)u,v}\rangle=D^2F(x)[u,v]= \langle{F''_e(x)u,v}\rangle_e 
=\langle{F''(e)F''_e(x)u,v}\rangle\]
yields $F''(x)=F''(e)F_e''(x)$, or 
\begin{equation}\label{eq:F''_e}
F_e''(x)=F''(e)^{-1}F''(x).
\end{equation}
Theorem~\ref{thm:scaling} implies that 
$F_e''(x)(K)=F''(e)^{-1}F''(x)(K)=F''(e)^{-1}(K^*)=K$, 
so that $F_e''(x)\in\Aut(K)$. Theorem~\ref{thm:scaling} also shows 
that, given any two points $u,v\in\Int(K)$, we can find a (unique)
point $z\in K$ such that $F''(z)u=F''(e)v\in K^*$. 
Therefore, $F_e''(x)(u)=v$, which shows that the set of linear operators 
$\{F_e''(x):x\in\Int(K)\}$ acts transitively on $K$. Hence, $K$ is a 
symmetric cone.

For the second assertion, note that if $s\in K^*_e=K$, then 
\[(F_e)_*(s):=\sup_{x\in K} \left\{ -\langle{x,s}\rangle_e-F(x)\right\}
=\sup_{x\in K}
\left\{-\langle{x,F''(e)s}\rangle-F(x)\right\}=F_*(F''(e)s).\]
For $x,z\in\Int(K)$, we thus have 
\[(F_e)_*(F''_e(z)x)=F_*(F''(e)F''_e(z)x)=F_*(F''(z)x)=F(x)-2F(z)-\nu,\]
where the second and last equalities follow from \eqref{eq:F''_e} and 
\eqref{eq:ss-2}, respectively. Consequently, $F$ is self--scaled
under $\langle{\cdot,\cdot}\rangle_e$. 
\end{proof}

\begin{remark}\label{rem:e}
Theorem~\ref{thm:K-symmetric} shows that from here on we may assume
without loss of generality that $K$ is a symmetric cone and that there
exists a (unique) point $e\in\Int(K)$ such that $F''(e)=\I$. 
\end{remark}

Together with Equation \eqref{eq:ss-2} this implies that 
\begin{equation}\label{eq:F_*=F} 
F_*(x)=F(x)-2F(e)-\nu=F(x)+const,
\end{equation}
and invoking \eqref{eq:ss-2} once more this implies the identity 
\begin{equation}\label{eq:sym-2}
F(F''(w)x)=F(x)-2F(w)+2F(e),\qquad\forall x,w\in\Int(K).
\end{equation}
Note that \eqref{eq:sym-2} is a criterion that involves only the
primal barrier $F$. Indeed, this identity allows one to characterise 
self--scaled barrier functions without invoking $F_*$, see
Lemma~\ref{lemma:charact-ss} below. Changing a barrier function 
by an additive constant is of no real consequence, as interior--point 
methods rely on gradient and Hessian information. Therefore, we can 
assume without loss of generality that $F(e)=0$. For the same reason,
Equation~\eqref{eq:F_*=F} makes it possible to think of $F$ and $F_*$
as the same function. Hence, we no longer need to distinguish 
between primal and dual quantities -- between $F$ and $F_*$, the 
primal and dual scaling points and so forth. 

We next prove a property of the Hessian $F''(w)$ which will become 
an essential tool for our classification of self--scaled barriers. 
For all $y\in\Int(K)$ let us define 
\[P(y):=F''(y)^{-1}.\]
\begin{lemma}\label{lemma:fund-eq}
For all $x,w\in\Int(K)$ it is true that 
\begin{equation}\label{eq:fund-formula}
P(P(w)x)=P(w)P(x)P(w). 
\end{equation} 
\end{lemma}

\begin{proof}
Let us define $z=P(w)x$. Equation \eqref{eq:sym-2} implies that for
any $h\in E$  we have $F(F''(w)(z+th))=F(z+th)-2F(w)+2F(e)$.
Expanding both sides of this equation and comparing the $t^2$ terms 
one gets \[D^2F(F''(w)z)[F''(w)h,F''(w)h]=D^2F(z)[h,h],\] or 
$\langle{F''(x)F''(w)h,F''(w)h }\rangle = \langle{F''(z)h,h}\rangle$. 
Thus, $F''(w)F''(x)F''(w)=F''(z)=F''(P(w)x)$, and 
\eqref{eq:fund-formula} follows.
\end{proof}
In the proof above, we need only the weaker condition 
$F(F''(w)x)=F(x)+c(w)$ where  $c(w)$ is a constant dependent on $w$.  
(However, Lemma~\ref{lemma:charact-ss} below shows that this is 
equivalent to \eqref{eq:sym-2}.) 
Equation \eqref{eq:fund-formula} is a symmetric version of a Formula
(3.2) from \cite{NT97}, see also Equation \eqref{eq:ff1stversion} above. 
In accordance with the established tradition in the theory of Jordan 
algebras we call \eqref{eq:fund-formula} the {\em fundamental formula}.

\begin{remark} 
We do not have a Jordan algebra connected to $F$ yet, but the 
fundamental formula leads one to suspect that there might be one. In 
spring 2000, inspired by the work of Petersson \cite{P95}, 
G\"uler~\cite{Guler00} proves that this is indeed the case. 
Subsequently, S.\ Schmieta~\cite{Schmieta00} independently discovers 
the same result, following essentially the same steps. Schmieta uses 
this fact as an essential tool to classify self--scaled barriers. As it 
turns out, the Jordan algebra connected to $F$ is  already
discovered by McCrimmon in his thesis \cite{McCrimmon65}, even without 
the assumption of convexity for $F$. His proof in turn is a 
generalisation of Koecher's ideas \cite{Koecher99} on
{\em $\omega$--domains}. Reading both works is instructive in delineating 
the role of convexity. 
\end{remark}

The following result provides an alternative definition of self--scaled 
barrier functions.

\begin{lemma}
\label{lemma:charact-ss}
Let $K$ be a regular, self--dual cone. A logarithmically homogeneous 
self--concordant barrier function $F$ on $\Int(K)$ is self-scaled if and only if 
\begin{align}
F''(w)x&\in\Int(K),\qquad\forall\,x,w\in\Int(K),\label{eq:sym-1'}\\ 
F(F''(w)x)&=F(x)+c(w),\qquad\forall\,x,w\in\Int(K),\label{eq:sym-2'}
\end{align}
where $c(w)$ is a constant that depends on $w$. 
\end{lemma}

\begin{proof} 
Since $K$ is self--dual, Equation \eqref{eq:sym-1'} is equivalent to  
Axiom \eqref{eq:ss-1}.  If $F$ is self--scaled, then 
Equation \eqref{eq:sym-2'} follows from \eqref{eq:sym-2}. 
Conversely, assume that $F$ satisfies Equation \eqref{eq:sym-2'}. 
Let $x,s\in\Int(K)$ be arbitrary points.  

We claim that there exists
a point $w\in\Int(K)$ such that $F''(w)x=s$.  Towards proving the claim, 
we consider the optimization 
problem $\min\{\langle{z^*,x}\rangle: \langle{z,s}\rangle=1\}$, where 
$z^*=-F'(z)$.  It is well known that the feasible region is bounded, 
see \cite{FK94}, Corollary I.1.6, p.~4. We have 
$F(x)+F_*(z^*)\ge-\nu-\nu\log\nu-\nu\log\langle{z^*,x}\rangle$ (see  
Proposition~\ref{prop2.1} $vii)$), and $F(z)+F_*(z^*)=-\nu$ 
(see Proposition~\ref{prop2.1} $vi)$), which imply
$F(x)-F(z)\ge-\nu\log\nu-\nu\log\langle{z^*,x}\rangle$.  These imply  
that the objective function of the optimization problem goes to infinity
as $z$ approaches the boundary of the feasible region, and thus    
the optimization problem has a minimizer $\hat{z}\in\Int(K)$ 
satisfying $F''(\hat{z})x=\lambda s$ for some scalar $\lambda$. Since 
$F''(\hat{z})x,s\in\Int(K)$, we have $\lambda>0$.  The point 
$w=\sqrt\lambda\hat{z}$ satisfies $F''(w)x=s$ (see 
Proposition~\ref{prop2.1} $v)$), and proves the claim.

Next, we claim that 
\begin{equation}
c(w)=-2F(w)+2F(e).\label{eq:c(w)}
\end{equation}
Let $u\in\Int(K)$ be a point satisfying $F''(u)w=e$. The fundamental 
formula \eqref{eq:fund-formula} is a consequence of 
\eqref{eq:sym-2'} and gives $F''(u)P(w)F''(u)=\I$, or 
equivalently $F''(w)=F''(u)^2$. From \eqref{eq:sym-2'}, we obtain  
\[ F(e) + c(w) = F(F''(w)e) = F(F''(u)^2e) = F(e) + 2c(u), \]
or $c(w)=2c(u)$. Equation \eqref{eq:sym-2'} also implies that  
\[ F(e) = F(F''(u)w) = F(w) + c(u) = F(w) +\frac{1}{2}c(w), \]
hence proving the claim. 

Using logarithmic homogeneity alone one can prove that
$F_*(w^*)=-\nu-F(w)$ where $w^*:=-F'(w)$ (see Proposition~\ref{prop2.1}
$vi)$).  Proposition~\ref{prop2.1} $ii)$ shows that the mapping 
$w\mapsto w^*$ is involutive, that is, $w^{**}=w$.  These imply 
$F_*(w)=F_*(w^{**})=-\nu-F(w^*)$. Since $F''(w)w=w^*$ by 
Proposition~\ref{prop2.1} $i)$, we have 
\[-\nu-F_*(w)=F(w^*)=F(F''(w)w)=F(w)-2F(w)+2F(e),\]
which is to say that $F_*(w)=F(w)-2F(e)-\nu$. This implies 
\[F_*(F''(w)x)=F(F''(w)x)-2F(e)-\nu=F(x)-2F(w)-\nu,\]
where the last equality follows from Equations \eqref{eq:sym-2'} and 
\eqref{eq:c(w)}. This concludes the proof.
\end{proof}
Note that together with Equation \eqref{eq:F_*=F},  
Lemma~\ref{lemma:charact-ss} 
implies that we can replace Axiom \eqref{eq:ss-2} of the original 
definition of a self--scaled barrier function by the requirement 
$F_*(F''(w)x)=F(x)+C(w)$ for some constant $C(w)$ which depends on 
$w$. This fact is already known, see \cite{Renegar99}.

\section{Group--Theoretic Aspects of Self--Scaledness}\label{sec:group}

In this section we explore the relationship between the Hessians of 
self--scaled barrier functions and the symmetry group of their domain
of definition. Though we present these results primarily for the 
purposes of later sections they are also of independent interest. 

The universal barrier function $U(x)$ defined in \eqref{eq:U(x)} 
plays an important role in the context of this section. The choice of 
the inner product $\langle\cdot,\cdot\rangle$ used in the definition 
of the characteristic function $\varphi_K(x)$ via \eqref{charfn} is 
irrelevant, since $\varphi_K$ changes only by an additive constant
under a change of  $\langle\cdot,\cdot\rangle$. G\"uler~\cite{Guler96} 
shows that the universal barrier function $U(x)$ is self--scaled, see 
Equation (13) and Theorem 4.4 in \cite{Guler96}. 
For all $x\in\Int(K)$ let \[Q(x):=U''(x)^{-1},\] and let $f\in\Int(K)$ 
be the point characterised by the equation \[Q(f)=\I.\]

\begin{remark}\label{rem:identity}
It follows from Theorem~\ref{thm:scaling} that $f$ is unique. The 
existence of such a point is also well known, see for example page 17 
of \cite{FK94}. 
\end{remark}

The point $f$ is the ``unit'' associated with the self--scaled 
barrier $U(x)$, see \cite{FK94} Proposition I.3.5, p.\ 14, and it 
is also the unit of the Jordan algebra associated with $U$. 

The following lemma is Theorem 3.17, pp.\ 205--206 in
\cite{Rothaus60}. We include a short proof of this result 
because these ideas play an important role in later sections. See 
also \cite{FK94}, Proposition I.4.3, p.\ 18 for a different approach 
to proving this result. 

\begin{lemma}\label{lemma:orthogonal} 
The orthogonal subgroup $\Orth\bigl(\Aut(K)\bigr)\subseteq\Aut(K)$ 
coincides with the stabiliser group at $f$, that is, 
\[\Orth\bigl(\Aut(K)\bigr)=\{H\in\Aut(K):Hf=f\}.\] 
\end{lemma}

\begin{proof}
If $A\in\Aut(K)$, then \[D^2U(Af)[Ah,Ah]=D^2U(f)[h,h]\] for every 
vector $h\in E$. That is, $A^*Q(Af)^{-1}A=\I$, or 
$Q(Af)=AA^*$, 
see for example \cite{Guler96}, Equation (11). This shows that $A$ is 
orthogonal if and only if $\I=Q(Af)$. The uniqueness of $f$ implies 
that this condition is equivalent to $Af=f$. 
\end{proof}

Next we note that the elements of $\Aut(K)$ have a unique 
{\em polar decomposition}, see \cite{Rothaus60}, Theorem 3.18, p.\ 206. 
For the sake of completeness we give a short proof.

\begin{lemma}\label{lemma:polar-decomposition}
Let $A\in\Aut(K)$. There exists a unique vector $u\in\Int(K)$ 
and a unique orthogonal cone automorphism 
$H\in\Orth\bigl(\Aut(K)\bigr)$ such that 
\[A=Q(u)H.\] 
\end{lemma}

\begin{proof}
By virtue of Theorem~\ref{thm:scaling}, there exists a unique point
$u\in\Int(K)$ such that $Q(u)f=Af$. Then $H:=Q(u)^{-1}A$ satisfies 
$Hf=Q(u)^{-1}Af=f$, which implies that $H$ is orthogonal by 
Lemma~\ref{lemma:orthogonal}. Since $H$ is orthogonal and $Q(u)$ is 
symmetric, $A=Q(u)H$ is indeed a polar decomposition of $A$. 

Suppose now that $A=Q_1H_1=Q_2H_2$ where $Q_i$ is
symmetric and $H_i$ is orthogonal, $i=1,2$. Then, 
$H:=Q_2^{-1}Q_1=H_2H_1^{-1}$ is orthogonal, and we have 
$\I=HH^*=Q_2^{-1}Q_1^2Q_2^{-1}$, or $Q_2^2=Q_1^2$. Since $Q_1$ and 
$Q_2$ are symmetric, we have $Q_1=Q_2$ and $H_1=H_2$. 
\end{proof}
  
The following result will play a key role in
Section~\ref{sec:classification} where we classify self--scaled
barriers. 

\begin{lemma}\label{lemma:Q=P}
The sets of inverse Hessians of $F$ and $U$ coincide, that is, 
\[\{P(v):v\in\Int(K)\}=\{Q(u):u\in\Int(K)\},\] 
and for all $x\in\Int(K)$ it is true that 
\begin{equation}\label{eq:PQ}
P(x)=Q(Q(x)^{1/2}e^{-1})=Q(x)^{1/2}Q(e)^{-1}Q(x)^{1/2},
\end{equation}
where $e^{-1}\in\Int(K)$ is characterised by the equation 
$Q(e^{-1})=Q(e)^{-1}$.
\end{lemma}

The point $e^{-1}$ is the inverse of $e$ in the Jordan algebra 
associated with $U(x)$. Note that Proposition~\ref{prop2.1} $iii)$ 
shows that $e^{-1}=-F'(e)$. 

\begin{proof}
If $v\in\Int(K)$, then Lemma~\ref{lemma:polar-decomposition} implies
that we can write $P(v)=Q(u)H$ 
for some $u\in\Int(K)$ and $H\in\Orth\bigl(\Aut(K)\bigr)$. 
By the uniqueness of the polar 
decompositions $P(v)=P(v)\I$ and $P(v)=Q(u)H$ it must be true that 
$P(v)=Q(u)$. Thus,  
\begin{equation*}
\{P(v):w\in\Int(K)\}\subseteq\{Q(u):u\in\Int(K)\}.
\end{equation*}
Conversely, let $u\in\Int(K)$. By Theorem~\ref{thm:scaling}, there 
exists a point $v\in\Int(K)$ such that $P(v)f=Q(u)f$. But this implies
that $Q(u)^{-1}P(v)f=f$. 
By virtue of Lemma~\ref{lemma:orthogonal} $H:=Q(u^{-1})P(v)$ is
therefore orthogonal. This means that $P(v)$ has the polar
decompositions $P(v)=P(v)\I$ and $P(v)=Q(u)H$. The uniqueness part of 
Lemma~\ref{lemma:orthogonal} then implies that $Q(u)=P(v)$ and
$H=\I$. This proves the first statement of the lemma.

Now let $x\in\Int(K)$ and define $u$ by $x=Q(u)f$, see 
Theorem~\ref{thm:scaling}. We have 
\begin{equation*}
Q(x)=Q\bigl(Q(u)f\bigr)=Q(u)Q(f)Q(u)=Q(u)^2,   
\end{equation*}
where the second equality follows from the fundamental Formula 
\eqref{eq:fund-formula}. In a similar vein, taking the first part of 
this lemma into account we obtain 
\begin{equation*}
P(x)=P\bigl(Q(u)f\bigr)=Q(u)P(f)Q(u).
\end{equation*}
These two equations imply that $Q(u)=Q(x)^{1/2}$ and 
\[P(x)=Q(x)^{1/2}P(f)Q(x)^{1/2}.\]
In particular, setting $x=e$ yields $\I=Q(e)^{1/2}P(f)Q(e)^{1/2}$,
that is, $P(f)Q(e)=\I$, 
and $P(f)=Q(e)^{-1}=Q(e^{-1})$. The lemma follows, since this 
implies that 
\[P(x)=Q(x)^{1/2}Q(e^{-1})Q(x)^{1/2}=Q(Q(x)^{1/2}e^{-1}).\]
\end{proof}

Although it does not have a direct bearing on later results, the
following proposition already shows that the self--scaled barrier
function $F$ is intimately connected to the universal barrier
function.   

\begin{proposition}
There exist constants $\alpha_1>0$ and $\alpha_2$ such that 
\[U(x)=\alpha_1\ln\determinant F''(x)+\alpha_2.\]
\end{proposition}

\begin{proof}
>From Equation \eqref{eq:PQ} we see that $\determinant P(x)
=\determinant Q(e)^{-1}\determinant Q(x)$, 
implying that $\determinant F''(x)=c_1\determinant U''(x)$ 
for some constant $c_1>0$. Theorem 4.4 in \cite{Guler96} shows 
that the function $u(x)=\ln\varphi_K(x)$ satisfies the equation 
$u(x)=c_2+\frac{1}{2}\ln\determinant u''(x)$ for some constant $c_2$. 
These facts combined with \eqref{eq:U(x)} imply the proposition. 
\end{proof} 

\section{Decomposition of Cones and Barrier Functions}\label{sec:decomp}

In this section, we prove two related results.  A cone is called 
{\em pointed} if it does not contain any whole lines.
First, we show that any pointed, convex cone decomposes into a 
direct sum of indecomposable or {\em irreducible} components in a unique
fashion. This result, which is of 
independent interest, is essentially a special case of Corollary 1 
in \cite{Gruber70}. Gruber's paper is the earliest mentioning of 
this result we could find, though it may have 
been derived several times independently. Gruber's original result 
addresses a more general affine setting which renders his proof more 
technically involved. Therefore, we include a simple and accessible 
proof. Second, we use this decomposition to write any self--scaled 
barrier function defined on the interior of a symmetric cone $K$ as a 
direct sum of self--scaled barriers defined on the irreducible 
components of $K$. 

Recall that the Minkowski sum of a set $\{A_i\}_{i=1}^k$ of subsets 
of $E$ is defined as 
\[A_1+\cdots+A_k:=\Bigl\{\sum_{i=1}^k x_i: x_i\in A_i\Bigr\}.\]
If all of the $A_i$ are linear subspaces $\{0\}\neq E_i\subseteq E$ 
which satisfy $E=E_1+\cdots+E_k$ and $E_i\cap(\sum_{j\neq i}E_j)=\{0\}$, 
then we say that the sum $E=E_1+\cdots+E_m$ is {\em direct} and write 
\[E=E_1\oplus E_2\oplus\cdots\oplus E_m.\]

\begin{definition}
Let $K\subseteq E$ be a pointed, convex cone. $K$ is called 
{\em decomposable} if there exist cones $\{K_i\}_{i=1}^m$, $m\ge2$, 
such that $K=K_1+\cdots+K_m$, where each $K_i$ ($i=1,\ldots,m$) 
lies in a linear subspace $E_i\subset E$, and where the spaces 
$\{E_i\}_{i=1}^m$ decompose $E$ into a direct sum 
$E=E_1\oplus E_2\oplus\cdots\oplus E_m$. Each $K_i$ is called a 
{\em direct summand} of $K$, and $K$ is called the {\em direct sum} 
of the $\{K_i\}$. We write 
\begin{equation}\label{eq:direct-sum} 
K=K_1\oplus K_2\oplus\cdots\oplus K_m
\end{equation}
to denote this relationship between $K$ and $\{K_i\}_{i=1}^m$. $K$ is
called {\em indecomposable} or {\em irreducible} if it cannot be
decomposed into a non-trivial direct sum. 
\end{definition}

Let us define $\hat{E}_i:=\oplus_{j\neq i}E_j$ and 
$\hat{K}_i:=\oplus_{j\neq i}K_j$. 
If $K$ is the direct sum \eqref{eq:direct-sum} then every 
$x\in K$ has a unique representation $x=x_1+\cdots+x_m$ with 
$x_i\in K_i\subseteq E_i$. Thus, 
$x_i=\pi_{E_i}x$, where $\pi_{E_i}$ is the projection of $E$ onto $E_i$ 
along $\hat{E}_i$. Also, since $0\in K_i$, we have 
$K_i =K_i+\sum_{j\neq i}\{0\} \subseteq \sum_{j=1}^m K_j=K$. 
Therefore, 
\[\pi_{E_i}K=K_i\subseteq K.\]
This implies that $K_i=\pi_{E_i}K$ is a convex cone.
Similarly, we have 
\[(\I-\pi_{E_i})K=\pi_{\hat{E}_i}K=\hat{K}_i\subseteq K.\]
We first prove a useful technical result:

\begin{lemma}\label{lemma:sum}
Let $K$ be a pointed, convex cone which decomposes into the direct sum 
\eqref{eq:direct-sum}. If $x\in K_i$ is a sum $x=x_1+\cdots+x_k$ 
of elements $x_j\in K$, then each $x_j\in K_i$. 
\end{lemma}

\begin{proof}  We have 
$0=\pi_{\hat{E}_i}x=\pi_{\hat{E}_i}x_1+\ldots+\pi_{\hat{E}_i}x_k$. 
Each term $\hat{x}_j:=\pi_{\hat{E}_i}x_j\in\hat{K}_i\subseteq K$, 
therefore we have $\hat{x}_j\in K$ and 
$-\hat{x}_j=\sum_{l\neq j}\hat{x}_l\in K$. Since $K$ contains no 
lines it must be true that $\hat{x}_j=0$, that is, 
$x_j=\pi_{E_i}x_j\in K_i$, ($j=1,\dots,k$).
\end{proof}

\begin{theorem}\label{uniqueness}
Let $K\subseteq E$ be a decomposable, pointed, convex cone.  The 
irreducible decompositions of $K$ are identical modulo indexing, 
that is, the set of cones $\{K_i\}_{i=1}^m$ is unique. Moreover, 
the subspaces $E_i$ corresponding to the non--zero cones $K_i$ are also 
unique. In particular, if $K$ is solid, then all the cones $K_i$ are 
non--zero and the subspaces $E_i$ are unique. 
\end{theorem}

\begin{proof}
Suppose that $K$ admits two irreducible decompositions 
\begin{align*}
K&=\bigoplus_{i=1}^m K_i\subseteq\bigoplus_{i=1}^m E_i\quad\text{ and}\\ 
K&=\bigoplus_{j=1}^q C_j\subseteq\bigoplus_{j=1}^q F_j.
\end{align*}
Note that each non--zero summand in either decomposition of $K$ must 
lie in $\Span(K)$ and that the subspace corresponding to each zero summand 
must be one--dimensional, for otherwise the summand would be decomposable. 
This implies that the number of zero summands in both decompositions is 
$\codim\bigl(\Span(K)\bigr)$. 
We may thus concentrate our efforts on $\Span(K)$, that is, 
we can assume that $K$ is solid and all the summands of both
decompositions of $K$ are non--zero. 
By \eqref{eq:direct-sum}, each $x\in C_j\subseteq K$ has a unique 
representation $x=x_1+\cdots+x_m$ where 
$x_i=\pi_{E_i}x\in K_i\subseteq K$. 
Also, Lemma~\ref{lemma:sum} implies that $x_i\in C_j$, and hence 
$x_i\in K_i\cap C_j$. Consequently, every 
$x\in C_j$ lies in the set $(K_1\cap C_j)+\cdots+(K_m\cap C_j)$. 
Conversely, we have $K_i\cap C_j\subseteq C_j$, implying 
that $(K_1\cap C_j)+\cdots+ (K_m\cap C_j)\subseteq C_j$.
Therefore, we have  
\[ C_j=(K_1\cap C_j)+\cdots+(K_m\cap C_j). \]
We have $K_i\cap C_j\subseteq E_i\cap F_j$, 
$F_j=(E_1\cap F_j)+\cdots+(E_m\cap F_j)$, and the intersection of 
any two distinct summands in the last sum is the trivial subspace 
$\{0\}$. The above decompositions of $F_j$ and $C_j$ are thus direct 
sums. Since $C_j$ is indecomposable, exactly one of the summands 
in the decomposition of $C_j$ is non--trivial. Thus, $C_j=K_i\cap C_j$, 
and hence $C_j\subseteq K_i$ for some $i$.  Arguing symmetrically, 
we also have $K_i\subseteq C_l$ for some $l$, implying that 
$C_j\subseteq C_l$. Therefore, $j=l$ for else 
$C_j\subseteq F_j\cap F_l=\{0\}$, contradicting our assumption above. 
This shows that $C_j=K_i$. The theorem is proved by repeating the above 
arguments for the cone 
$\hat{K}_i=\oplus_{k\neq i}K_k=\oplus_{l\neq j}C_l$. 
\end{proof}

Next, we show that self--scaled barrier functions have irreducible
decompositions as well. Let $F$ be defined on $\Int(K)$ where 
$K$ is a symmetric cone with irreducible decomposition
\eqref{eq:direct-sum}. For $i=1,\dots,m$, let $F_i$ be a function 
defined on  $\relint(K_i)$, the relative interior of $K_i$ in
$E$. If $F(x)=\sum_{i=1}^m F_i(x_i)$ 
for every $x=\sum_{i=1}^m x_i\in\oplus_{i=1}^m\relint(K_i)=\Int(K)$, 
then we say that $F$ is the direct sum of the $F_i$ and write 
\begin{equation}\label{eq:dirsum-F}
F=\bigoplus_{i=1}^m F_i.
\end{equation}

\begin{theorem}\label{thm:decomp-F}
Let $K$ be a symmetric cone with irreducible decomposition 
\eqref{eq:direct-sum}. The irreducible components $K_i$ are 
symmetric cones. Let $F(x)$ be a self--scaled barrier for $K$. 
There exist self--scaled barrier functions $F_i$ for the 
cones $K_i$ such that 
\begin{equation*} 
F=F_1\oplus\dots\oplus F_m.
\end{equation*}
The functions $F_i$ are unique up to additive constants. 
\end{theorem}

\begin{proof}
Recall that the universal barrier function $U(x)$ is 
given in \eqref{eq:U(x)}. Changing the inner product used in 
the definition of the characteristic function $\varphi_K(x)$ changes 
$U(x)$ only by an additive constant, hence we may assume that 
$\langle{x,y}\rangle=\sum_{i=1}^m\langle{x_i,y_i}\rangle_{E_i}$ for the
purposes of this definition. Here, $\langle\cdot,\cdot\rangle_{E_i}$ 
is an inner product defined on $E_i$ chosen so that $U_i''(f_i)=id_{E_i}$ 
for some elements $f_i\in\relint(K_i)$ where $U_i$ denotes the universal 
barrier function defined on $\relint(K_i)$. Then we have $Q(f)=\I$ for 
$f=\oplus_{i=1}^m f_i\in\Int(K)$, in full consistency with our 
previous notation. Moreover, $K$ is self--dual under
$\langle\cdot,\cdot\rangle$, since $K^*=Q(f)K=K$. Hence, 
we may choose the vector $e\in\Int(K)$ specified in Remark~\ref{rem:e}
as the unique element in $\Int(K)$ such that $F''(e)=\I$ under
$\langle\cdot,\cdot\rangle$, see Remark~\ref{rem:identity}. 

Thus, $U$ can be written as the direct sum 
$U(x)=\oplus_{i=1}^m U_i(x_i)$ and $Q(x)$ has block structure
corresponding to the subspaces $E_i$, 
$Q(x)=\oplus_{i=1}^m Q_i(x_i)$ where $Q_i(x_i)=U_i''(x_i)^{-1}$. 
Consequently, Equation \eqref{eq:PQ} implies that $P(x)$ also has the 
same block structure, 
\begin{equation}\label{eq:blockstructure}
P(x)=P_1(x_1)\oplus\dots\oplus P_m(x_m),
\end{equation} 
where $P_i(x)=Q_i\bigl(Q_i(x_i)^{1/2}e_i^{-1}\bigr)\in\Aut(K_i)$ 
with $e_i^{-1}=\pi_{E_i}\bigl(-F'(e)\bigr)$, $\pi_{E_i}$ being the
projection at the beginning of this section. 

So far we know that $P(x)$ has block structure corresponding to the 
direct sum $E=\oplus_{i=1}^m E_i$, but it is not a priori clear 
that $P_i(x_i)^{-1}$ is the Hessian of a function defined on
$\relint(K_i)$. Let the spaces $\hat{E}_i$ be 
defined as earlier in this section, and let us consider the vector fields 
$v_i:x\mapsto\pi_{E_i}F'(x)$, defined on $\Int(K)$ and taking values
in $E_i$ for all $i$, $i=1,\ldots,m$. We claim that $v_i$ depends only on
$x_i=\pi_{E_i}x$. In fact, for any two vectors $x,y\in\Int(K)$ 
such that $x_i=y_i$ we have 
\begin{align*}
v_i(y)&=\pi_{E_i}F'(y)
=\pi_{E_i}\Bigl[F'(x)+\int_{0}^{1}F''\bigl(ty+(1-t)x\bigr)[y-x]dt\Bigr]\\
&=v_i(x)+\int_{0}^{1}P_i\bigl(\pi_{E_i}\bigl[ty+(1-t)x\bigr]\bigr)^{-1}
\pi_{E_i}[y-x]dt=v_i(x)+\int_{0}^{1}0\,dt,
\end{align*}
which shows our claim. Hence, the quotient vector fields 
\begin{align*}
\hat{v}_i:\Int(K)/\hat{E}_i&\rightarrow E_i,\\
x\mod\hat{E}_i&\mapsto\pi_{E_i}F'(x)
\end{align*} 
are well defined and can be identified with vector fields $\hat{v}_i$ 
defined on the cones $\relint(K_i)$. The direct sum of these vector
fields amounts to the gradient field 
\begin{equation}\label{eq:vectorfields}
F'=\hat{v}_i\oplus\dots\oplus\hat{v}_m:\Int(K)\rightarrow E.
\end{equation}
$F'$ being conservative, the $\hat{v}_i$ must be conservative too, 
implying that these are the gradient vector fields of some functions 
$F_i$ defined on $\relint(K_i)$ which are uniquely determined up to 
additive constants. We may choose these constants so that 
$F=\oplus_{i=1}^m F_i$. Clearly, we have $F_i''(x_i)=P_i(x)^{-1}$ 
for any $x\in\Int(K)$. 

Using Equation~\eqref{eq:blockstructure}, it is straightforward to
check that the $F_i$ are self-concordant, see \cite{NN94}. Applying 
Proposition~\ref{prop2.1} $i)$ and $iv)$ to $F$, using 
Equation~\eqref{eq:vectorfields} and considering variations of 
$x\in\Int(K)$ only in the part $x_i=\pi_{E_i}x$, we obtain that 
$\langle x_i,-F_i'(x_i)\rangle=\nu_i$ for some number $\nu_i>0$. 
Moreover, applying Proposition~\ref{prop2.1} $v)$ to 
$F$ and using Equation~\eqref{eq:vectorfields} we get 
$F_i'(\tau x_i)=\tau^{-1}F_i'(x_i)$ for all $\tau>0$. Hence, 
\begin{gather*}
\begin{split}
F_i\bigl(\tau x_i\bigr)&=F_i(x_i)+\int_1^\tau
\frac{d}{d\xi}F_i\bigl(\xi x_i\bigr)d\xi=
F_i(x_i)+\int_1^\tau\bigl\langle x_i,F_i'(\xi x_i)\bigr\rangle d\xi\\
&=F_i(x_i)-\int_1^\tau\xi^{-1}\bigl\langle x_i,-F_i'(x_i)\bigr\rangle d\xi
=F_i(x_i)-\nu_i\int_1^\tau\xi^{-1}d\xi\\
&=F_i(x_i)-\nu_i\ln\tau.
\end{split}
\end{gather*}
This shows that the functions $F_i$ are
$\nu_i$\nobreakdash-logarithmically homogeneous. It is a well--known 
fact that any logarithmically homogeneous self--concordant function is
also a barrier function, see for example \cite{NN94} or \cite{Renegar99}. 
It remains to show that the functions $F_i$ are self--scaled. As previously 
noted, Condition~\eqref{eq:sym-1'} is satisfied, since 
$P_i(x_i)\in\Aut(K_i)$ for all $i$. Finally, Condition 
\eqref{eq:sym-2'} holds for $F_i$ because we can apply this condition 
to $F$, choosing $w=w_i\oplus\bigl(\oplus_{j\neq i}e_i\bigr)$ and 
$x=x_i\oplus\bigl(\oplus_{j\neq i}e_i\bigr)$ and using the block 
structures of $F$ and $F''$. 

Note that the irreducible components $K_i$ of $K$ must be symmetric
cones, since the $F_i$ are self--scaled barriers defined on
$\relint(K_i)$. The symmetry of the $K_i$ can also 
be directly derived from the block structure of
$Q(x)=\oplus_{i=1}^mQ_i(x_i)$ and the fact that the set of cone
automorphisms $\bigl\{Q(x):x\in\Int(K)\bigr\}$ acts transitively on
$K$. 
\end{proof}

The decomposition Theorem~\ref{thm:decomp-F} shows that for the purposes
of
classifying self--scaled barriers we may concentrate our efforts on 
irreducible cones. 

\section{Classification of Self--Scaled Barriers}
\label{sec:classification}

In this section, we give a complete classification of self--scaled 
barrier functions on the symmetric cone $K$. 

The definition of a self--scaled barrier function $F$ requires that 
$F$ changes only by an additive constant under the action of symmetric
cone automorphisms $\{P(u):u\in\Int(K)\}$, see Equation
\eqref{eq:sym-2}. However, it is not a priori known how $F$ behaves
under the action of an arbitrary element of $\Aut(K)$. Note that this
is in marked contrast to the case where $F$ is the universal barrier
function $U$, which is known to change only by an additive constant
under the action of any element of the symmetry group of $K$. This 
explains in a sense the main difficulty one faces when proving the 
results below. 

The next result is key in resolving this difficulty and is just a 
slight reformulation of the conjecture raised in \cite{Hauser99},
according to which self-scaled barriers on irreducible symmetric cones
are isotropic.  Let us denote by $\Aut(K)_0$ the connected component of
the identity in $\Aut(K)$. 

\begin{theorem}\label{thm:isotropic}
Let $K$ be a symmetric cone. If $H\in\Aut(K)_0$ is 
orthogonal, then $F(Hx)=F(x)$ for all $x\in\Int(K)$. 
\end{theorem}

\begin{proof}
Koecher~\cite{Koecher99} proves that if $K$ is a symmetric cone, then 
$\Aut(K)_0$ is generated by $\{Q(u):u\in{\mathcal V}\}$ where 
${\mathcal V}$ is a neighbourhood of the identity $f$, see \cite{Koecher99},
Theorem 4.9 (b), pp.\ 88--89. Koecher's proof
exploits the fact that all 
{\em derivations} of the Jordan algebra associated with $U(x)$ are 
{\em inner}. An accessible proof for the case where $K$ is irreducible
is given in \cite{FK94}, Lemma VI.1.2, pp.\ 101--102, based on certain 
non--trivial results from the theory of Jordan algebras. If 
$H\in\Aut(K)_0$ is orthogonal, it follows from Koecher's result that 
\[H=\prod_1^l Q(u_i)=\prod_1^l P(v_i),\] for some $u_i, v_i\in\Int(K)$, 
$(i=1,\dots,l)$. Here the second equality follows from
Lemma~\ref{lemma:Q=P}. 
Therefore, it follows from Equation \eqref{eq:sym-2} that 
\[F(Hx)=F(\prod_1^l P(v_i) x)=F(x)-2\sum_1^l F(v_i).\]
Since $Hf=f$, setting $x=f$ above yields $\sum_1^l F(v_i)=0$ and
settles the claim of the theorem. 
\end{proof}

The group $\Aut(K)_0$ already acts transitively on $K$, see \cite{FK94}, 
page 5. Thus, the above result is significant.

\begin{remark}\label{rem:lim}
Hauser's approach ~\cite{Hauser00} to solving the isotropy conjecture 
for the cone of positive semidefinite symmetric matrices is based on 
similar ideas. Hauser essentially rediscovers Koecher's Theorem 4.9
(b) in this special case, see \cite{Hauser00} Proposition 3.3 and 
Corollary 4.3. He then uses Proposition 3.3 in conjunction with the 
fundamental formula as the key mechanism in the proof of the
conjecture. Y.\ Lim \cite{LimPC} generalised this approach to
arbitrary irreducible 
symmetric cones while refereeing the paper \cite{Hauser00}, thus 
completing the classification of self--scaled barriers. However, 
since Lim was an anonymous referee, his result was not publicly 
announced. 
\end{remark}

Now, assume that $K$ is
irreducible and let $U$ be the universal barrier function for $K$. 
Let $k$ be the rank of the Jordan algebra associated with $U(x)$. Let
$x$ be an arbitrary point in $\Int(K)$. Then there exists an
orthogonal frame $\{e_1,\ldots,e_k\}$ such that 
$x=\sum_{i=1}^k\lambda_i e_i$, $\lambda_i>0$, $i=1,\ldots,k$, 
see \cite{FK94}, Theorem~III.1.2, pp.\ 44--45. 
By $C$ let us denote the cone generated by this frame, that is, 
\[C=\bigl\{\sum_{i=1}^k\lambda_i e_i:\lambda_i\ge 0\bigr\}.\]
Note that $C$ is a direct sum of the half--lines 
$\{\lambda_i e_i:\lambda_i\ge0\}$ and thus a symmetric cone. 

\begin{lemma}\label{lemma:restriction}
If $F$ is a self--scaled barrier function defined on the interior of
the irreducible symmetric cone $K$, and if $C$ is as defined above, then 
\[ F(\sum_{i=1}^k\alpha_ie_i) = -\frac{\nu}{k}\sum_{i=1}^k\log\alpha_i,
\qquad\alpha_i>0\quad(i=1,\ldots,k). \]
\end{lemma}

\begin{proof}
Let $\sigma$ be any permutation of $\{1,\ldots,k\}$.  Theorem IV.2.5 
in \cite{FK94} implies that there exists an orthogonal automorphism 
$H\in Aut(K)_0$ such that $He_i=e_{\sigma(i)}$ ($i=1,\ldots,k$). 
Using Theorem~\ref{thm:isotropic}, we then obtain 
\[ F(\sum_{i=1}^k\alpha_ie_{\sigma(i)})
=F(\sum_{i=1}^k\alpha_ie_i),\qquad\forall\alpha_i>0\quad(i=1,\ldots,k). \]
Define $g(\alpha_1,\ldots,\alpha_k):=F(\sum_{i=1}^k\alpha_ie_i)$. 
Note that $g$ is a symmetric function.
Consider a point $e+\sum_{i=1}^k\beta_ie_i=\sum_{i=1}^k\alpha_i e_i
\in\relint(C)$, with arbitrary $\beta_i\ge 0$ and $\alpha_i:=1+\beta_i$. 
Applying Theorem 5.1 in \cite{NT97} repeatedly, we obtain  
\[ F(e+\sum_{i=1}^k\beta_ie_i)-F(e) = \sum_{i=1}^k (F(e+\beta_ie_i)-F(e)). \]
Using the symmetry of $g$ and $F(e)=0$, the above equation translates into 
\[ g(\alpha_1,\ldots,\alpha_k) 
  = \sum_{i=1}^k g(\alpha_i,1,\ldots,1),\qquad\forall\alpha_i\ge1\quad(i=1,\ldots,k). \]
If $\alpha_1=\ldots=\alpha_k=\alpha$ above, we have 
$g(\alpha,\ldots,\alpha)=F(\alpha e)=F(e)-\nu\log\alpha=-\nu\log\alpha$. 
This gives $g(\alpha,1,\ldots,1)=-\frac{\nu}{k}\log\alpha$ for all 
$\alpha\ge1$.  Consequently, we have 
\begin{equation}
\label{eq:g}
  g(\alpha_1,\ldots,\alpha_k) 
= -\frac{\nu}{k}\sum_{i=1}^k\log\alpha_i,\qquad
\forall\alpha_i\ge1\quad(i=1,\ldots,k).
\end{equation}
Now, if $\alpha_i>0$ ($i=1,\ldots,k$) are arbitrary, choose $t>0$ such 
that $\hat{\alpha}_i=\alpha_i/t\ge1$.  Since $F$ is logarithmically 
homogeneous, we have 
$g(\alpha_1,\ldots,\alpha_k)=g(\hat{\alpha}_1,\ldots,\hat{\alpha}_k)-\nu\log t$ 
by the logarithmic homogeneity of $F$.  A simple calculation shows that 
\eqref{eq:g} holds true for all $\alpha_i>0$.  The lemma is proved.
\end{proof}

The following theorem classifies self--scaled barrier functions for 
irreducible symmetric cones. 

\begin{theorem}\label{thm:classify-1}
Let $K$ be an irreducible symmetric cone, and let $F$ be a self--scaled 
barrier function defined on $\Int(K)$. Then there exist constants 
$\alpha>0$ and $\beta$ such that 
\[F(x)=\alpha U(x)+\beta,\]
where $U(x)$ is the universal barrier function on $\Int{K}$. 
\end{theorem}

\begin{proof}
Lemma~\ref{lemma:restriction} describes the restriction of $F$ on $\relint(C)$. 
Since the universal barrier function is also self--scaled, the same 
considerations apply to $U(x)-U(e)$. Thus, the functions $F$ and $U$
are homothetic transformations of each other on $\relint(C)$, that is, 
there exist $\alpha>0$, $\beta$ such that 
\begin{equation}\label{eq:homothetic}
F(x)=\alpha U(x)+\beta.
\end{equation}
Let $y\in\Int(K)$ be an arbitrary point with the spectral decomposition 
$y=\sum_{i=1}^k\nu_i d_i$.  Corollary IV.2.7 in \cite{FK94} implies 
that there exists $A\in O(\Aut(K)_0)$ such that $Ae_i=d_i$ 
($i=1,\ldots,k$).  We have $y=Ax$ where $x=\sum_{i=1}^k\nu_i e_i\in\relint(C)$.
Theorem~\ref{thm:isotropic} gives $F(y)=F(x)$ and $U(y)=U(x)$, and hence 
the Identity \eqref{eq:homothetic} extends to all of $\Int(K)$. 
\end{proof}

We are now ready to give the final classification theorem for
self--scaled barrier functions on arbitrary symmetric cones. This
theorem shows that all self--scaled barrier functions are related 
to the standard logarithmic or the universal barrier via homothetic
transformations. 

\begin{theorem}\label{thm:main}
Let $F$ be a self--scaled barrier function for a symmetric cone $K$ 
with irreducible decomposition \eqref{eq:direct-sum}. 
Then there exist constants $c_0$ and $c_1\ge1,\dots,c_m\ge 1$ such that 
\begin{equation*}
F=c_0-\bigoplus_{i=1}^m c_i\ln\determinant_{K_i},
\end{equation*}
see Equation \eqref{eq:dirsum-F}.  Here $\determinant_{K_i}x_i$ 
denotes the determinant of $x_i\in\relint(K_i)$ in the Jordan 
algebraic sense, see \cite{FK94}, Chapter 2. Conversely, any function 
of this form is a self--scaled barrier for $K$. 
\end{theorem}

\begin{proof}
Theorems \ref{thm:decomp-F} and \ref{thm:classify-1} imply that there 
exist constants $d_0$ and $d_1>0,\dots,d_m>0$ such that 
\[F(x)=d_0+d_1 u_1(x_1)+\ldots+d_m u_m(x_m),\]
where $u_i(x_i)=\ln\varphi_{K_i}(x_i)$. It is known that 
$u_i(x_i)=const-n_i/r_i\ln\determinant x_i$, where $r_i$ is the rank of 
the Jordan algebra associated with $u_i(x)$, and $n_i$ is the 
dimension of the cone $K_i$, see \cite{FK94}, Proposition~III.4.3, 
p.\ 53. Finally, Theorem~4.1 in \cite{GT98} implies that the 
function $-\alpha\ln\determinant x_i$ is self--concordant if and only if 
$\alpha\ge1$. 
\end{proof}

\section{Acknowledgements}

R.\ Hauser conducted his research on this paper partially while
writing his Ph.D.\ at the School of Operations Research and Industrial
Engineering of Cornell University, and partially as a postdoctoral
fellow at the University of Cambridge. He wishes to thank
Mike Todd, Jim Renegar and Arieh Iserles for numerous valuable 
discussions and for their continuous support. This research was 
supported by a ``bourse pour chercheur d\'ebutant'' from the Swiss 
National Science Foundation and the Swiss Academy of Technical
Sciences, and by the Engineering and Physical Sciences Research 
Council of the UK, EPSRC grant GR/M30975. 

The research of O.\ G\"uler on this paper was conducted at the
Industrial Engineering Department at Bilkent University, Ankara, 
Turkey, while on a sabbatical leave from the Department of Mathematics
and Statistics, University of Maryland, Baltimore County, Baltimore, 
Maryland. The author thanks the Industrial Engineering Department 
at Bilkent University for providing him a congenial atmosphere and
excellent 
working conditions. This research was partially supported by the National 
Science Foundation under grant DMS--0075722.

\end{document}